\documentclass[11pt]{amsart}
 \usepackage{amssymb,amsfonts,amscd,amsbsy}
\usepackage[mathscr]{eucal}
\usepackage{url}

\theoremstyle{plain}
\newtheorem{thm}{Theorem}[section]
\newtheorem{lem}[thm]{Lemma}
\newtheorem{prop}[thm]{Proposition}
\newtheorem{cor}[thm]{Corollary}
\theoremstyle{definition}

\numberwithin{equation}{section}

\begin{document}
\title[$M_2$-rank differences for overpartitions]{$M_2$-rank differences for overpartitions}
\author{Jeremy Lovejoy and Robert Osburn}

\address{CNRS, LIAFA, Universit{\'e} Denis Diderot - Paris 7, Case 7014, 75205 Paris Cedex 13, FRANCE}

\address{School of Mathematical Sciences, University College Dublin, Belfield, Dublin 4, Ireland}

\email{lovejoy@liafa.jussieu.fr}

\email{robert.osburn@ucd.ie}



\dedicatory{In memory of Oliver Atkin}

\date{\today}

\begin{abstract}
This is the third and final installment in our series of papers applying the method of Atkin and Swinnerton-Dyer to deduce 
formulas for rank differences.  The study of rank differences was initiated by Atkin and Swinnerton-Dyer in their
proof of Dyson's conjectures concerning Ramanujan's congruences for the partition
function. Since then, other types of rank differences for statistics associated to partitions
have been investigated. In this paper, we prove explicit formulas for $M_2$-rank differences for
overpartitions. Additionally, we express a third order mock theta function in terms of rank differences.
\end{abstract}

\maketitle
\section{Introduction}
The rank of a partition $\lambda$ is defined to be the largest part $\ell(\lambda)$ minus the number of parts $n(\lambda)$.  Let $N(s,m,n)$ denote the number of partitions
of $n$ with rank congruent to $s$ modulo $m$.  Responding to a conjecture of Dyson \cite{dyson}, Atkin and Swinnerton-Dyer proved elegant formulas, in terms of modular functions and generalized Lambert series, for the generating functions for $N(r,\ell,\ell n + d) - N(s,\ell,\ell n+d)$ when $\ell = 5$ and $7$.  
For example, they found \cite[Theorem 4]{asd}:

\begin{equation} \label{asdexample1}
\sum_{n \geq 0} \bigl(N(1,5,5n+2) - N(2,5,5n+2)\bigr)q^n = \frac{(q^5;q^5)_{\infty}}{(q^2,q^3;q^5)_{\infty}}
\end{equation}

\noindent and

\begin{equation} \label{asdexample2}
\begin{aligned}
\sum_{n \geq 0} \bigl(N(0,5,5n+3) &- N(1,5,5n+3) \bigr)q^n  \\ 
&= \frac{q}{(q^5;q^5)_{\infty}}\sum_{n \in \mathbb{Z}} \frac{(-1)^nq^{15n(n+1)/2}}{1-q^{5n+2}} + \frac{(q,q^4,q^5;q^5)_{\infty}}{(q^2,q^3;q^5)_{\infty}^2}.
\end{aligned}
\end{equation}

Here we have employed the standard basic hypergeometric series
notation \cite{Ga-Ra1},
$$
(a_1,a_2,\dots,a_j;q)_n = \prod_{k=0}^{n-1}(1-a_1q^k)(1-a_2q^k)
\cdots (1-a_jq^k).
$$
We follow the custom of dropping the ``$;q$" unless the base is something other than $q$.  

The rank of a partition studied by Atkin and Swinnerton-Dyer is now understood to be a special case of a more general rank which is defined on overpartition pairs \cite{Br-Lo1}.  Recall that an overpartition of $n$ is a partition of $n$ where we may overline the first occurrence of a part while an overpartition pair $(\lambda,\mu)$ of $n$ is a pair of overpartitions where the sum of all of the parts is $n$. The rank of an overpartition pair $(\lambda,\mu)$ is

\begin{equation} \label{pairrank}
\ell((\lambda,\mu)) - n(\lambda) - \overline{n}(\mu) -
\chi((\lambda,\mu)),
\end{equation}

\noindent where $\overline{n}(\cdot)$ is the number of overlined parts only and $\chi((\lambda,\mu))$ is defined to be $1$ if the largest part of $(\lambda,\mu)$ occurs only non-overlined and only in $\mu$, and $0$ otherwise.  

When $\mu$ is empty and $\lambda$ has no overlined parts, \eqref{pairrank} becomes the rank of a partition.  In addition to this rank, three other special cases of \eqref{pairrank} have turned out to be of particular interest: the rank of an overpartition, the $M_2$-rank of a partition without repeated odd parts, and the $M_2$-rank of an overpartition.  For more on these three ranks and their generating functions, see \cite{bg, Br-Lo.5, Br-Lo2, Br-Lo-Os.5, Br-Lo-Os1, Br-Zw1, Dewar, Lo1, Lo2, loveoz, Lo-Os1}.
In \cite{loveoz} and \cite{Lo-Os1}, we applied the method of Atkin and Swinnerton-Dyer to find formulas like \eqref{asdexample1} and \eqref{asdexample2} for certain rank differences for overpartitions and $M_2$-rank differences for partitions without repeated odd parts.  Here we complete the picture by doing the same for $M_2$-rank differences for overpartitions.     

This rank arises by replacing parts $m$ in the first component of an overpartition pair $(\lambda,\mu)$ by $2m$ and parts $m$ in the second component by $2m-1$.  From \eqref{pairrank}, the $M_2$-rank of the resulting overpartition $\lambda$ is
$$
\text{$M_2$-rank}(\lambda) := \bigg \lceil \frac{\ell(\lambda)}{2}
\bigg \rceil - n(\lambda) + n(\lambda_o) - \chi(\lambda),
$$
where $\lambda_o$ is the subpartition consisting of the odd non-overlined parts and 
$\chi(\lambda) = 1$ if the largest part of $\lambda$ is odd
and non-overlined and $\chi(\lambda) = 0$ otherwise.  For example,
the $M_2$-rank of the overpartition
$5+\overline{4}+4+\overline{3}+1+1$ is $3-6+3-1 = -1$.

Let $\overline{N}_{2}(s, \ell, n)$ denote the number of overpartitions of $n$ whose $M_{2}$-rank is congruent to $s$ modulo $\ell$. 
Using the notation
\begin{equation} \label{rst}
R_{st}(d) = \sum_{n \geq 0} \left(\overline{N}_2(s,\ell,\ell n + d) -
\overline{N}_2(t,\ell,\ell n + d)\right ) q^n,
\end{equation}
where the prime $\ell$ will always be clear, our main results are
summarized in Theorems \ref{main3} and \ref{main} below.

\begin{thm}\label{main3}
For $\ell=3$, we have

\begin{equation} \label{r_01(0)}
R_{01}(0) =  -1 + \frac{(-q)_{\infty}(q^3;q^3)_{\infty}^2}{(q)_{\infty}(-q^3;q^3)_{\infty}^2},
\end{equation}

\begin{equation} \label{r_01(1)}
R_{01}(1) = \frac{2(q^3;q^3)_{\infty}(q^6;q^6)_{\infty}}{(q)_{\infty}},
\end{equation}

\begin{equation}
R_{01}(2) = \frac{4(q^6;q^6)_{\infty}^4}{(q^2;q^2)_{\infty}(q^3;q^3)_{\infty}^2} + \frac{6q(-q^3;q^3)_{\infty}}{(q^3;q^3)_{\infty}}\sum_{n \in \mathbb{Z}} \frac{(-1)^nq^{3n^2+6n}}{1-q^{6n+2}}.
\end{equation}
\end{thm}

\begin{thm} \label{main} For $\ell=5$, we  have

\begin{equation} \label{r_{12}(0)}
R_{12}(0)= \frac{10q^2(q^{10};q^{10})_{\infty}^4(q,q^2,q^8,q^9;q^{10})_{\infty}}{(q)_{\infty}^3(q;q^2)_{\infty}} + \frac{2q(q^5;q^5)_{\infty}}{(q;q^2)_{\infty}^5(q^3,q^4,q^6,q^7;q^{10})_{\infty}},
\end{equation}

\begin{equation} \label{r_{12}(1)}
\begin{aligned}
R_{12}(1) & = \frac{-6q^3(-q^5;q^5)_{\infty}}{(q^5;q^5)_{\infty}} \sum_{n \in \mathbb{Z}} \frac{(-1)^nq^{5n^2+10n}}{1-q^{10n+4}} + \frac{4q(q^{10};q^{10})_{\infty}^2}{(q;q^2)_{\infty}^5(q^4,q^6;q^{10})_{\infty}(q^5;q^5)_{\infty}} \\ 
& + \frac{20q^3(q^{10};q^{10})_{\infty}^7(q,q^9;q^{10})_{\infty}^2(q^2,q^8;q^{10})_{\infty}^3}{(q)_{\infty}^4(q^5;q^5)_{\infty}^2},
\end{aligned}
\end{equation}

\begin{equation} \label{r_{12}(2)}
R_{12}(2)= \frac{10q(q^{10};q^{10})_{\infty}^3}{(q)_{\infty}^2(q^3,q^7;q^{10})_{\infty}^3(q,q^2,q^8,q^9;q^{10})_{\infty}},
\end{equation}

\begin{equation} \label{r_{12}(3)}
\begin{aligned}
R_{12}(3) &= \frac{10(q^{10};q^{10})_{\infty}^3}{(q)_{\infty}^2(q^3,q^4,q^6,q^7;q^{10})_{\infty}(q,q^9;q^{10})_{\infty}^3} \\
&- \frac{8(q^5;q^5)_{\infty}}{(q;q^2)_{\infty}^5(q^2,q^3,q^7,q^8;q^{10})_{\infty}},
\end{aligned}
\end{equation}

\begin{equation} \label{r_{12}(4)}
\begin{aligned}
R_{12}(4) &=  \frac{-2q(-q^5;q^5)_{\infty}}{(q^5;q^5)_{\infty}} \sum_{n \in \mathbb{Z}} \frac{(-1)^nq^{5n^2+10n}}{1-q^{10n+2}} + \frac{4(q^2;q^2)_{\infty}}{(q;q^2)_{\infty}^4(q^2,q^3,q^7,q^8;q^{10})_{\infty}^3},
\end{aligned}
\end{equation}

\begin{equation} \label{r_{02}(0)}
\begin{aligned}
R_{02}(0) &= -1 + \frac{(q^5;q^5)_{\infty}^6}{(q;q^2)_{\infty}^6(q^3,q^4,q^6,q^7;q^{10})_{\infty}(q^{10};q^{10})_{\infty}^5} \\
& + \frac{q^2(q)_{\infty}(q^{10};q^{10})_{\infty}}{(q^3,q^4,q^6,q^7;q^{10})_{\infty}^3(q^{5};q^{5})_{\infty}} + \frac{4q(q^5;q^5)_{\infty}^2(q,q^9;q^{10})_{\infty}}{(q;q^2)_{\infty}^6(q^4,q^6,q^{10};q^{10})_{\infty}} \\
& - \frac{10q^2(q^{10};q^{10})_{\infty}^3}{(q)_{\infty}^2(q^2,q^8;q^{10})_{\infty}(q^3,q^7;q^{10})_{\infty}^4},
\end{aligned}
\end{equation}

\begin{equation} \label{r_{02}(1)}
\begin{aligned}
R_{02}(1) &= \frac{2q^3(-q^5;q^5)_{\infty}}{(q^5;q^5)_{\infty}} \sum_{n \in \mathbb{Z}} \frac{(-1)^nq^{5n^2+10n}}{1-q^{10n+4}} + \frac{2(q^{5};q^{5})_{\infty}^4}{(q;q^2)_{\infty}^6(q^4,q^6;q^{10})_{\infty}(q^{10};q^{10})_{\infty}^3} \\ 
& + \frac{2q^2(q)_{\infty}(q^{10};q^{10})_{\infty}^3}{(q^3,q^7;q^{10})_{\infty}^2(q^4,q^6;q^{10})_{\infty}^3(q^5;q^5)_{\infty}^3},
\end{aligned} 
\end{equation}

\begin{equation} \label{r_{02}(2)}
R_{02}(2)= \frac{4(q^3,q^7;q^{10})_{\infty}(q^5;q^5)_{\infty}^2}{(q;q^2)_{\infty}^6(q^4,q^6,q^{10};q^{10})_{\infty}} - \frac{10q(q^{10};q^{10})_{\infty}^3}{(q)_{\infty}^2(q^3,q^7;q^{10})_{\infty}^3(q,q^2,q^8,q^9;q^{10})_{\infty}},
\end{equation}

\begin{equation} \label{r_{02}(3)}
R_{02}(3)= \frac{4(q^5;q^5)_{\infty}}{(q;q^2)_{\infty}^5(q^2,q^3,q^7,q^8;q^{10})_{\infty}} ,
\end{equation}

\begin{equation} \label{r_{02}(4)}
\begin{aligned}
R_{02}(4) &= \frac{4q(-q^5;q^5)_{\infty}}{(q^5;q^5)_{\infty}} \sum_{n \in \mathbb{Z}} \frac{(-1)^nq^{5n^2+10n}}{1-q^{10n+2}} -\frac{10(q^{10};q^{10})_{\infty}^3(q^5;q^5)_{\infty}(q^4,q^6;q^{10})_{\infty}}{(q,q^9;q^{10})_{\infty}(q)_{\infty}^3(q^3,q^7;q^{10})_{\infty}} \\
&- \frac{2(q^{10};q^{10})_{\infty}(-q^2,-q^3;q^5)_{\infty}}{q(-q,-q^4;q^5)_{\infty}(q,q^2,q^8,q^9;q^{10})_{\infty}} + \frac{2(q^5;q^5)_{\infty}^2(q^2,q^8;q^{10})_{\infty}^5(q^{10};q^{10})_{\infty}^2}{q(q;q^2)_{\infty}^4(q,q^9;q^{10})_{\infty}^4(q^2;q^2)_{\infty}^3}.
\end{aligned}
\end{equation}

\end{thm}

The method of Atkin and Swinnerton-Dyer may be generally described as
regarding groups of identities as equalities between polynomials
of degree $\ell-1$ in $q$ whose coefficients are power series in
$q^{\ell}$. Specifically, we first consider the expression

\begin{equation} \label{term}
\sum_{n=0}^{\infty} \Bigl\{ \overline{N}_{2}(s,\ell,n) -
\overline{N}_{2}(t,\ell,n) \Bigr\} q^{n}
\frac{(q)_{\infty}}{2(-q)_{\infty}}.
\end{equation}

\noindent By (\ref{gen1}), (\ref{s}), and (\ref{final}), we write
(\ref{term}) as a polynomial in $q$ whose coefficients are power
series in $q^{\ell}$.  We then alternatively express (\ref{term})
in the same manner using the formulas in Theorem \ref{main} 
and equation \eqref{lem6eq1} or \eqref{lem6eq2}.
Finally, we use the theory of modular forms to show that these
two resulting polynomials are the same for each pair of values of
$s$ and $t$.  

Some comments are in order here.  First, if the
number of overpartitions of $n$ with $M_{2}$-rank $m$ is denoted by
$\overline{N}_{2}(m,n)$, then $\overline{N}_{2}(m,n) = \overline{N}_{2}(-m,n)$ (see \eqref{gen}). 
Hence the values of $s$ and $t$ considered in Theorems \ref{main3} and
\ref{main} are sufficient to find any rank difference generating
function $R_{st}(d)$.  Second, the formulas in Theorems \ref{main3} and \ref{main} are somewhat more
complicated than the ones in \cite{asd}, \cite{loveoz}, and \cite{Lo-Os1}.  When there are exactly two infinite products, we have verified using Euler's algorithm \cite[p. 98, Ex. 2]{An1} that they cannot be reduced to one product. However, in \eqref{r_{02}(0)} and \eqref{r_{02}(4)} we cannot rule out the possibility of a simpler expression.  

Finally, the formulas for $R_{01}(0)$ and $R_{02}(1)$ when $\ell = 3$ match those for the classical rank
differences for overpartitions \cite[Eqs. (1.1) and (1.2)]{loveoz}.  In other words, letting $\overline{N}(s,m,n)$ denote the number of overpartitions whose rank is $s$ modulo $m$, we have 
$$
\overline{N}_2(0,3,3n+d) - \overline{N}_2(1,3,3n+d) = \overline{N}(0,3,3n+d) - \overline{N}(1,3,3n+d)
$$
for $n \geq 0$ and $d=0$ or $1$.  When $d=2$ it turns out that the generating function for the difference of the rank differences is 
proportional to the third order mock theta function
$$
\omega(q):=\sum_{n=0}^{\infty} \frac{q^{2n(n+1)}}{(q; q^2)_{n+1}^2}.
$$




\begin{cor} \label{omega} 
We have
\begin{equation*} \label{mock}
\begin{aligned}
6\omega(q) &=\sum_{n \geq 0} \bigl( \overline{N}_2(0,3,3n+2) - \overline{N}_2(1,3,3n+2) \bigr)q^n \\
& - \sum_{n \geq 0} \bigl( \overline{N}(0,3,3n+2) - \overline{N}(1,3,3n+2) \bigr)q^n.
\end{aligned}
\end{equation*}
\end{cor}

This is not the first time that mock theta functions have appeared in relation to rank differences.   Andrews and Garvan \cite[Section 4]{An-Ga1} and Hickerson \cite[Section 5]{Hi1} have already shown that certain fifth and seventh order mock theta functions can be expressed in terms of rank differences of Atkin and Swinnerton-Dyer.  Some tenth order mock theta functions are also rank differences.  Specifically, using identities for the tenth order mock theta functions $\phi(q)$ and $\psi(q)$ on pages $533$--$534$ of \cite{Ch1}, combined with identities $(1.9)$, $(1.11)$, and $(1.14)$ of \cite{loveoz}, we have

\begin{equation*}
2\phi(q) = \sum_{n \geq 0} \bigl( \overline{N}(0,5,5n+1) - \overline{N}(2,5,5n+1) \bigr) q^n
\end{equation*}

\noindent and

\begin{equation*}
2\psi(q) = \sum_{n \geq 0} \bigl( \overline{N}(0,5,5n+4) + \overline{N}(1,5,5n+4) - 2\overline{N}(2,5,5n+4) \bigr)q^{n+1}.
\end{equation*} 
In general, the generalized Lambert series which arise in the study of rank differences are known to be building blocks of mock theta functions \cite{Zw1}.


The paper is organized as follows.  In Section 2 we collect some
basic definitions, notations and generating functions.  
In Section 3 we prove two key $q$-series
identities relating generalized Lambert series to infinite
products, and in Section 4 we give the proofs of Theorems
\ref{main3} and \ref{main}. In Section 5, we prove Corollary \ref{omega}.


\section{Preliminaries}
We begin by introducing some notation and definitions, essentially
following \cite{asd}.  With $y=q^{\ell}$, let

$$
r_{s}(d):=\sum_{n=0}^{\infty} \overline{N}_{2}(s,\ell, \ell n+d) y^n
$$
and
$$
r_{st}(d):=r_{s}(d) - r_{t}(d).
$$
Thus we have
$$
\sum_{n=0}^{\infty} \overline{N}_{2}(s,\ell,n) q^n =
\sum_{d=0}^{\ell-1} r_{s}(d) q^d.
$$

To abbreviate the sums occurring in Theorems \ref{main3} and \ref{main}, we define
$$
\Sigma(z,\zeta, q):= \sum_{n \in \mathbb{Z}} \frac{(-1)^n
\zeta^{2n} q^{n^2 + 2n}}{1-z^{2} q^{2n}}.
$$
Henceforth we assume that $a$ is not a multiple of $\ell$.  We
write
$$
\Sigma(a,b):=\Sigma(y^a, y^b, y^{\ell})=\sum_{n \in \mathbb{Z}}
\frac{(-1)^n y^{2bn + \ell n(n+2)}}{1-y^{2\ell n + 2a}}
$$
and
$$
\Sigma(0,b) :=  \sideset{}{'} \sum_{n \in \mathbb{Z}} \frac{(-1)^n
y^{2bn + \ell n(n+2)}}{1-y^{2\ell n}},
$$
where the prime means that the term corresponding to $n=0$ is
omitted.

To abbreviate the products occurring in Theorems \ref{main3} and \ref{main}, we
define

$$
P(z,q):=\prod_{r=1}^{\infty} (1-zq^{r-1})(1-z^{-1}q^r)
$$

\noindent and

$$
P(0):=\prod_{r=1}^{\infty} (1-y^{2\ell r}).
$$

\noindent We also have the relations

\begin{equation} \label{p1}
P(z^{-1}q, q)=P(z,q)
\end{equation}

\noindent and

\begin{equation} \label{p2}
P(zq, q)=-z^{-1}P(z,q).
\end{equation}

In \cite{Lo2}, it is shown that the two-variable generating
function for $\overline{N}_{2}(m,n)$ is

\begin{equation} \label{gen}
\sum_{n=0}^{\infty} \overline{N}_{2}(m,n) q^n = \frac{2(-q)_{\infty}}
{(q)_{\infty}} \sum_{n=1}^{\infty} (-1)^{n+1} q^{n^2 + 2|m|n}
\frac{1-q^{2n}}{1+q^{2n}}.
\end{equation}
From this we may easily deduce that the generating function for
$\overline{N}_{2}(s,m,n)$ is

\begin{equation} \label{gen1}
\sum_{n=0}^{\infty} \overline{N}_{2}(s,m,n) q^n =
\frac{2(-q)_{\infty}}{(q)_{\infty}} \sideset{}{'} \sum_{n \in
\mathbb{Z}} \frac{(-1)^n q^{n^2 + 2n}(q^{2sn} +
q^{2(m-s)n})}{(1+q^{2n})(1 - q^{2mn})}.
\end{equation}
Hence it will be beneficial to consider sums of the form

\begin{equation} \label{s}
\overline{S}_{2}(b):= \sideset{}{'}\sum_{n \in \mathbb{Z}}
\frac{(-1)^n q^{n^2+2bn}}{1-q^{2\ell n}}.
\end{equation}
We will require the relation

\begin{equation} \label{rels}
\overline{S}_{2}(b)=-\overline{S}_{2}(\ell-b),
\end{equation}
which follows from the substitution $n \to -n$ in (\ref{s}). We
shall also exploit the fact that the functions
$\overline{S}_{2}(\ell)$ are essentially infinite products.

\begin{lem}\label{Sofq}
We have
$$
\overline{S}_{2}(\ell) = \frac{-(q)_{\infty}}{2(-q)_{\infty}} +
\frac{1}{2}.
$$
\end{lem}

\begin{proof}
Use the relation \eqref{rels} to compute $-2\overline{S}_{2}(\ell)$ then apply the case $z = -1$ of Jacobi's triple product identity,

\begin{equation} \label{jtp}
\sum_{n \in \mathbb{Z}} z^nq^{n^2} = (-zq,-q/z,q^2;q^2)_{\infty}.
\end{equation}

\end{proof}

\section{Two lemmas}
The proofs of Theorems \ref{main3} and \ref{main} will follow from
identities which relate the sums $\displaystyle \Sigma(a,b)$ to
the products $P(z, q)$. The key steps are the two lemmas below.  The first is equation (5.4) of \cite{Br-Lo-Os1}
 
\begin{lem} \label{jack} We have

\begin{equation} \label{chan}
\begin{aligned}
& \sum_{n=-\infty}^{\infty} (-1)^n q^{n^2 + 2n} \Bigl[ \frac{\zeta^{-2n}}{1-z^{2}{\zeta^{-2}}q^{2n}} +
\frac{\zeta^{2n+4}}{1-z^{2}{\zeta^{2}}q^{2n}} \Bigr] \\
& = \frac{-2 (\zeta^4, q^{2}\zeta^{-4}; q^2)_{\infty}(-q)_{\infty}^2}{(-\zeta^2, -q\zeta^{-2})_{\infty} (\zeta^{-2}, q^2 \zeta^2 ; q^2)_{\infty}}
\sum_{n=-\infty}^{\infty} (-1)^n \frac{q^{n^2 + 2n}}{1-z^2 q^{2n}} \\
& + \frac{(-z^2, -qz^{-2})_{\infty} (\zeta^4, q^2 \zeta^{-4}, \zeta^2, q^2 \zeta^{-2}; q^2)_{\infty} (q^2; q^2)_{\infty}^2}
{(-\zeta^2, -q\zeta^{-2})_{\infty} (z^2 \zeta^2, q^2 z^{-2} \zeta^{-2}, z^2 \zeta^{-2}, q^{2} z^{-2} \zeta^2, z^2, z^{-2} q^2; q^2)_{\infty}}.
\end{aligned}
\end{equation}

\end{lem}


We now specialize Lemma \ref{chan} to the case $\zeta=y^a$, $z=y^b$, and $q=y^{\ell}$:

\begin{equation} \label{lem1}
\begin{aligned}
y^{4a} \Sigma(a+b, a) + \Sigma(b-a, -a) +& \frac{P(y^{4a}, y^{2\ell})P(-1, y^{\ell})}{P(-y^{2a}, y^{\ell}) P(y^{-2a}, y^{2\ell})} \Sigma(b,0) \\
& - \frac{P(-y^{2b}, y^{\ell})P(y^{4a}, y^{2\ell})P(y^{2a},
y^{2\ell})P(0)^2}{P(y^{2b+2a}, y^{2\ell})P(y^{2b-2a}, y^{2\ell})P(-y^{2a}, y^{\ell})P(y^{2b}, y^{2\ell})} = 0.
\end{aligned}
\end{equation}

\noindent We now define

$$
\begin{aligned}
g(z,q) & :=-\frac{P(z^4, q^2)P(-1,q)}{P(-z^2,q)P(z^{2}q^2,q^2)}\Sigma(z,1,q) - z^4 \Sigma(z^2, z, q) \\
& -  \sideset{}{'} \sum_{n=-\infty}^{\infty} \frac{(-1)^n z^{-2n} q^{n(n+2)}}{1-q^{2n}}
\end{aligned}
$$

\noindent and

\begin{equation} \label{g}
g(a):=g(y^a, y^{\ell})=-\frac{P(y^{4a}, y^{2\ell})P(-1, y^{\ell})}{P(-y^{2a}, y^{\ell})P(y^{-2a}, y^{2\ell})} \Sigma(a,0) - y^{4a} \Sigma(2a, a) - \Sigma(0, -a).
\end{equation}

\noindent The second key lemma is the following.

\begin{lem} \label{us2} We have

\begin{equation} \label{part1}
2g(z,q) - g(z^2, q) + \frac{1}{2} =\frac{P(z^6, q^2)^2 (q^2; q^2)_{\infty}^2}{P(z^2, q^2)^2 P(z^8, q^2)} - \frac{P(z^2, q)^2 P(z^4, q) (q)_{\infty}^2}{P(-z^2, q)^2 P(-z^4, q)P(-1,q)}
\end{equation}

\noindent and

\begin{equation} \label{part2}
g(z,q) + g(z^{-1}q, q)=0.
\end{equation}

\end{lem}

\begin{proof}

We first require a short computation involving $\Sigma(z, \zeta,
q)$. Note that

\begin{equation} \label{sigma}
\begin{aligned}
z^2 \Sigma(z,\zeta, q) &+ q\zeta^2 \Sigma(zq, \zeta, q) \\
& =\sum_{n=-\infty}^{\infty} (-1)^n \frac{z^2 \zeta^{2n} q^{n(n+2)}}{1 - z^2 q^{2n}} +
\sum_{n=-\infty}^{\infty} (-1)^{n} \frac{\zeta^{2n+2} q^{n(n+2) + 1}}{1-z^2 q^{2n+2}} \\
& = -\sum_{n=-\infty}^{\infty} (-1)^n \zeta^{2n} q^{n^2}
\end{aligned}
\end{equation}

\noindent upon writing $n-1$ for $n$ in the second sum of the
first equation. Taking $\zeta=1$ yields

\begin{equation} \label{step}
\begin{aligned}
z^2 \Sigma(z,1,q) + q\Sigma(zq, 1, q) & = -\sum_{n=-\infty}^{\infty} (-1)^n q^{n^2}. \\
\end{aligned}
\end{equation}

\noindent Now write $g(z,q)$ in the form

$$
g(z,q)=f_{1}(z) - f_{2}(z) - f_{3}(z)
$$

\noindent where

$$
f_{1}(z):= -\frac{P(z^4, q^2)P(-1,q)}{P(-z^2,q)P(z^2 q^2,q^2)}\Sigma(z,1,q),
$$

$$
f_{2}(z):= z^4 \Sigma(z^2, z, q),
$$

\noindent and

$$
f_{3}(z):= \sideset{}{'} \sum_{n=-\infty}^{\infty} \frac{(-1)^n
z^{-2n} q^{n(n+2)}}{1-q^{2n}}.
$$

\noindent By (\ref{p1}), (\ref{p2}), and (\ref{step}),

\begin{equation} \label{f1}
\begin{aligned}
f_{1}(zq) - f_{1}(z) & = \sum_{n=-\infty}^{\infty} (-1)^n q^{n^2} \frac{P(z^4, q^2) P(-1, q)}{P(-z^2,q) P(z^2,q^2)}.  \\
\end{aligned}
\end{equation}

\noindent A similar argument as in (\ref{sigma}) yields

\begin{equation} \label{f2}
f_{2}(zq) - f_{2}(z) = \sum_{n=-\infty}^{\infty} (-1)^n z^{2n} q^{n^2}
\end{equation}

\noindent and

\begin{equation} \label{f3}
f_{3}(zq) - f_{3}(z) = -1 + \sum_{n=-\infty}^{\infty} (-1)^n z^{-2n}
q^{n^2}.
\end{equation}

\noindent Adding (\ref{f2}) and (\ref{f3}), then subtracting from (\ref{f1}) gives

\begin{equation} \label{constant}
g(z,q)-g(zq, q) =-1.
\end{equation}

\noindent Here we have used the identity

\begin{equation} \label{hidden1}
\begin{aligned}
& \sum_{n=-\infty}^{\infty} (-1)^n q^{n^2} \frac{P(z^4, q^2) P(-1, q)}{P(-z^2,q) P(z^2,q^2)} = \sum_{n=-\infty}^{\infty} (-1)^n z^{2n}
q^{n^2} + \sum_{n=-\infty}^{\infty} (-1)^n z^{-2n} q^{n^2} \\
\end{aligned}
\end{equation}

\noindent which follows from the triple product identity \eqref{jtp} after writing $-n$ for $n$ in the second sum. If we now define

$$
f(z):= 2g(z,q) - g(z^2, q) + \frac{1}{2} - \frac{P(z^6, q^2)^2 (q^2; q^2)_{\infty}^2}{P(z^2, q^2)^2 P(z^8, q^2)} + \frac{P(z^2, q)^2 P(z^4, q) (q)_{\infty}^2}{P(-z^2, q)^2 P(-z^4, q)P(-1,q)}
$$

\noindent then from (\ref{p1}), (\ref{p2}), and (\ref{constant}),
one can verify that

\begin{equation} \label{functional}
f(zq) - f(z) =0.
\end{equation}

\noindent  Now, it follows from a routine complex analytic argument
similar to the proof of Lemma 4.2 in \cite{loveoz} (see also Lemma
2 in \cite{asd}) that $f(z)=0$. This proves (\ref{part1}).

To prove (\ref{part2}), it suffices to
show, after (\ref{constant}),

\begin{equation} \label{gees}
g(z^{-1}, q) + g(z,q)=-1.
\end{equation}

\noindent Note that

\begin{equation} \label{short}
\begin{aligned}
z^2\Sigma(z,1, q) + z^{-2} \Sigma(z^{-1}, 1, q) & =
z^2 \sum_{n=-\infty}^{\infty} (-1)^n \frac{q^{n(n+2)}}{1 - z^2 q^{2n}} -
\sum_{n=-\infty}^{\infty} (-1)^n \frac{q^{n^2}}{1-z^2 q^{2n}} \\
& = -\sum_{n=-\infty}^{\infty} (-1)^n q^{n^2}
\end{aligned}
\end{equation}

\noindent where we have written $-n$ for $n$ in the second sum in the first equation. Thus, by
(\ref{p1}), (\ref{p2}), and (\ref{short}), we have

\begin{equation} \label{f1z}
f_{1}(z) + f_{1}(z^{-1}) = -\sum_{n=-\infty}^{\infty} (-1)^n q^{n^2} \frac{P(z^4,
q^2)P(-1,q)}{P(-z^2, q)P(z^2, q^2)}.
\end{equation}

\noindent Again, a similar argument as in (\ref{short}) gives

\begin{equation} \label{f2z}
f_{2}(z) + f_{2}(z^{-1}) = -\sum_{n=-\infty}^{\infty} (-1)^n z^{2n} q^{n^2}
\end{equation}

\noindent and

\begin{equation} \label{f3z}
f_{3}(z) + f_{3}(z^{-1}) = 1 - \sum_{n=-\infty}^{\infty} (-1)^n z^{-2n} q^{n^2}.
\end{equation}

\noindent Adding (\ref{f2z}) and (\ref{f3z}), then subtracting from (\ref{f1z})
yields (\ref{gees}). Here we have again used (\ref{hidden1}).

\end{proof}

Letting $z=y^a$ and $q=y^{\ell}$ in Lemma \ref{us2}, we get

\begin{equation} \label{g1}
2g(a) - g(2a) + \frac{1}{2} = \frac{P(y^{6a}, y^{2\ell})^2 P(0)^2}{P(y^{2a}, y^{2\ell})^2 P(y^{8a}, y^{2\ell})} - \frac{P(y^{2a}, y^{\ell})^2 P(y^{4a}, y^{\ell}) (y^{\ell}; y^{\ell})_{\infty}^2}{P(-y^{2a}, y^{\ell})^2 P(-y^{4a}, y^{\ell}) P(-1, y^{\ell})}
\end{equation}

\noindent and

\begin{equation} \label{g2}
g(a) + g(\ell-a)=0.
\end{equation}

\noindent These two identities will be of key importance in the
proofs of Theorems \ref{main3} and \ref{main}.

\section{Proofs of Theorems \ref{main3} and \ref{main}}

We now compute the sums $\overline{S}_{2}(\ell-m)$. The reason for
this choice is two-fold. First, we would like to obtain as simple
an expression as possible in the final formulation (\ref{final}). Secondly, to prove Theorem \ref{main3}, we only need to compute $\overline{S}_{2}(1)$
whereas to prove Theorem \ref{main}, we need $\overline{S}_{2}(1)$ and $\overline{S}_{2}(3)$. The former yields $\overline{S}_{2}(2)$ while the latter in turn yields $\overline{S}_{2}(4)$ via (\ref{rels}). For $\ell=3$, we can choose $m=1$ and for $\ell=5$, $m=1$ and $m=2$ respectively.
As this point, we follow the idea of Section 6 in \cite{asd}. Namely, we write

\begin{equation} \label{n}
n=\ell r + m + b,
\end{equation}

\noindent where $-\infty < r < \infty$. The idea is to simplify
the exponent of $q$ in $\overline{S}_{2}(\ell-m)$. Thus

$$
2\ell n-2mn + n^2 = \ell^{2} r(r+2) + 2b\ell r + (b+m)(b-m + 2\ell).
$$

\noindent We now substitute (\ref{n}) into (\ref{s}) and let $b$
take the values $0$, $\pm a$, and $\pm m$. Here $a$ runs through
$1$, $2$, \dotso, $\frac{\ell-1}{2}$ where the value $a \equiv \pm
m \bmod \ell$ is omitted. As in \cite{asd}, we use the notation
$\displaystyle \sideset{}{''} \sum_{a}$ to denote the sum over
these values of $a$. We thus obtain

$$
\begin{aligned}
\overline{S}_{2}(\ell-m) & = \sideset{}{'} \sum_{n=-\infty}^{\infty} (-1)^n \frac{q^{2(\ell-m)n + n^2}}{1-y^{2n}} \\
& = \sum_{b} {\sideset{}{'} \sum_{r=-\infty}^{\infty}}
(-1)^{r+m+b} q^{(b+m)(b-m+2\ell)} \frac{y^{\ell r(r+2) + 2br}}{1-
y^{2\ell r + 2m + 2b}},
\end{aligned}
$$

\noindent where $b$ takes values $0$, $\pm a$, and $\pm m$ and the term corresponding
to $r=0$ and $b=-m$ is omitted. Thus

\begin{equation} \label{s(b)}
\begin{aligned}
\overline{S}_{2}(\ell-m) & = (-1)^m q^{m(2\ell-m)} \Sigma(m,0) + \Sigma(0, -m) + y^{4m} \Sigma(2m, m) \\
& + \sideset{}{''} \sum_{a} (-1)^{m+a} q^{(a+m)(a-m + 2\ell)} \Bigl
\{ \Sigma(m+a, a) + y^{-4a} \Sigma(m-a, -a) \Bigr \}.
\end{aligned}
\end{equation}

\noindent Here the first three terms arise from taking $b=0$, $-m$, and $m$ respectively.
We now can use (\ref{lem1}) to simplify this expression. By taking $b=m$ and dividing by
$y^{4a}$ in (\ref{lem1}), the sum of the two terms inside the curly brackets becomes

$$
\begin{aligned}
-y^{-4a} \frac{P(y^{4a}, y^{2\ell})P(-1, y^{\ell})}{P(-y^{2a}, y^{\ell})P(y^{-2a}, y^{2\ell})}
& \Sigma(m,0) \\
& + y^{-4a} \frac{P(-y^{2m}, y^{\ell})P(y^{4a}, y^{2\ell})P(y^{2a}, y^{2\ell})
P(0)^2}{P(y^{2a+2m}, y^{2\ell})P(y^{2m-2a}, y^{2\ell})P(-y^{2a}, y^{\ell})P(y^{2m}, y^{2\ell})}.
\end{aligned}
$$

\noindent Similarly, upon taking $a=m$ in (\ref{g}), then the sum of the second and third terms in (\ref{s(b)}) is

$$
-\frac{P(y^{4m}, y^{2\ell})P(-1, y^{\ell})}{P(-y^{2m}, y^{\ell})P(y^{-2m}, y^{2\ell})} \Sigma(m,0)
- g(m).
$$

In total, we have

\begin{equation} \label{final}
\begin{aligned}
\overline{S}_{2}(\ell-m) & = -g(m) \\
& + \sideset{}{''} \sum_{a} \Biggl\{ (-1)^{m+a} q^{(a+m)(a-m+2\ell)} y^{-4a} \\
& \times \frac{P(-y^{2m}, y^{\ell})P(y^{4a}, y^{2\ell})P(y^{2a}, y^{2\ell})P(0)^2}{P(y^{2a + 2m}, y^{2\ell})P(y^{2m-2a}, y^{2\ell})P(-y^{2a}, y^{\ell})P(y^{2m}, y^{2\ell})} \Biggr \} \\
& + \Sigma(m,0) \Biggl\{ (-1)^m q^{m(2\ell-m)} - \frac{P(y^{4m}, y^{2\ell})P(-1, y^{\ell})}{P(-y^{2m}, y^{\ell})P(y^{-2m}, y^{2\ell})} \\
& - \sideset{}{''} \sum_{a} (-1)^{m+a} q^{(a+m)(a-m+2\ell)} y^{-4a}
\frac{P(y^{4a}, y^{2\ell})P(-1, y^{\ell})}{P(-y^{2a}, y^{\ell})}P(y^{-2a}, y^{2\ell}) \Biggr\}.
\end{aligned}
\end{equation}

\noindent We can simplify some of the terms appearing in
(\ref{final}) as we are interested in certain values of $\ell$,
$m$, and $a$. To this end, we prove the following result. Let
$\{ \quad \}$ denote the coefficient of $\Sigma(m,0)$ in (\ref{final}).

\begin{prop} \label{brackets} If $\ell=3$ and $m=1$, then

$$
\{ \quad \} = -q^5 \frac{(q)_{\infty} (-q^{9};
q^{9})_{\infty}}{(-q)_{\infty} (q^{9}; q^{9})_{\infty}}.
$$

\noindent If $\ell=5$, $m=2$, and $a=1$, then

$$
\{ \quad \}=q^{16} \frac{(q)_{\infty} (-q^{25}; q^{25})_{\infty}}{(-q)_{\infty} (q^{25}; q^{25})_{\infty}}.$$

\noindent If $\ell=5$, $m=1$, $a=2$, then

$$
\{ \quad \}=-q^9  \frac{(q)_{\infty} (-q^{25};
q^{25})_{\infty}}{(-q)_{\infty} (q^{25}; q^{25})_{\infty}}.
$$

\end{prop}

\begin{proof}
This is a straightforward application of the identities
\begin{equation} \label{lem6eq1}
\frac{(q)_{\infty}}{(-q)_{\infty}} = \frac{(q^{9};
q^{9})_{\infty}}{(-q^{9}; q^{9})_{\infty}}
 - 2q(q^{3}, q^{15}, q^{18}; q^{18})_{\infty}
\end{equation}
and
\begin{equation} \label{lem6eq2}
\frac{(q)_{\infty}}{(-q)_{\infty}} = \frac{(q^{25};
q^{25})_{\infty}}{(-q^{25}; q^{25})_{\infty}}
 - 2q(q^{15}, q^{35}, q^{50}; q^{50})_{\infty} + 2q^4 (q^{5}, q^{45}, q^{50}; q^{50})_{\infty}.
\end{equation}
These are Lemma 3.1 in \cite{loveoz}.
\end{proof}

We are now in a position to prove Theorems \ref{main3} and
\ref{main}.  We begin with Theorem \ref{main3}.

\begin{proof}
By (\ref{gen1}), (\ref{s}), and (\ref{rels}), we have

\begin{equation} \label{gen3too}
\sum_{n=0}^{\infty} \Bigl\{ \overline{N}_{2}(0,3,n) -
\overline{N}_{2}(1,3,n) \Bigr\} q^{n}
\frac{(q)_{\infty}}{2(-q)_{\infty}}=3\overline{S}_{2}(1) +
\overline{S}_{2}(3).
\end{equation}

\noindent  By (\ref{p1}), (\ref{p2}), (\ref{final}), and
Proposition \ref{brackets} we have

\begin{equation} \label{s1too}
\overline{S}_{2}(1)= g(1) + q^2 y\Sigma(1,0) \frac{(q)_{\infty}
(-q^{9}; q^{9})_{\infty}} {(-q)_{\infty} (q^{9}; q^{9})_{\infty}}.
\end{equation}

By Lemma \ref{Sofq} we have
\begin{equation}
\overline{S}_{2}(3) = \frac{-(q)_{\infty}}{2(-q)_{\infty}} + \frac{1}{2}.
\end{equation}

We have that

$$
 3g(1) + 3q^2 y\Sigma(1,0) \frac{(q)_{\infty}
(-q^{9}; q^{9})_{\infty}} {(-q)_{\infty} (q^{9}; q^{9})_{\infty}}
- \frac{(q)_{\infty}}{2(-q)_{\infty}} + \frac{1}{2}  = \Bigl \{
r_{01}(0) q^0 + r_{01}(1)q + r_{01}(2)q^2  \Bigr \}
\frac{(q)_{\infty}}{2(-q)_{\infty}}.
$$

\noindent We now multiply the right hand side of the above
expression using \eqref{lem6eq1} and the $R_{01}(d)$ from Theorem
\ref{main} (recall that $r_{01}(d)$ is just $R_{01}(d)$ with $q$
replaced by $q^{3}$).  We then equate coefficients of powers of
$q$ and verify the resulting identities.  For $q^1$ and $q^2$ the 
resulting equation follows easily upon cancelling factors in infinite products. 
For $q^0$ we obtain
$$
3g(1) + \frac{1}{2} = 
\frac{(-q^3;q^3)_{\infty} (q^9;q^9)_{\infty}^3}{2(q^3;q^3)_{\infty}(-q^9;q^9)_{\infty}^3}
- 4y
\frac{(q^{18};q^{18})_{\infty}^4(q^3, q^{15},q^{18}; q^{18})_{\infty}}{(q^6;q^{6})_{\infty}(q^9;q^9)_{\infty}^2}.
$$
Appealing to \eqref{g1} and \eqref{g2} and then replacing $q$ by $q^{1/3}$, this becomes
$$
\frac{P(q^2,q^3)^2P(q^4,q^3)(q^3;q^3)_{\infty}^2}{P(-q^2,q^3)^2P(-q^4,q^3)P(-1,q^3)} + \frac{(-q;q)_{\infty} (q^3;q^3)_{\infty}^3}{2(q;q)_{\infty}(-q^3;q^3)_{\infty}^3}
= 4y
\frac{(q^{6};q^{6})_{\infty}^4(q, q^{5},q^{6}; q^{6})_{\infty}}{(q^2;q^{2})_{\infty}(q^3;q^3)_{\infty}^2}.
$$
After making a common denominator on the left and simplifying, this equation may be verified using the case $(z,\zeta,t,q) = (-q^2,q^2,-1,q^3)$ of the addition theorem \cite[Eq. (3.7)]{asd},
$$
P^2(z,q)P(\zeta t,q)P(\zeta/t,q) - P^2(\zeta,q)P(zt,q)P(z/t,q) +
\zeta/t P^2(t,q)P(z\zeta,q)P(z/\zeta,q) = 0.
$$  
This completes the proof of Theorem \ref{main3}.
\end{proof}

We now turn to Theorem \ref{main}.

\begin{proof}
We begin with the rank differences $R_{12}(d)$.  By (\ref{gen1}),
(\ref{s}), and (\ref{rels}), we have

\begin{equation} \label{gen3}
\sum_{n=0}^{\infty} \Bigl\{ \overline{N}_{2}(1,5,n) -
\overline{N}_{2}(2,5,n) \Bigr\} q^{n}
\frac{(q)_{\infty}}{2(-q)_{\infty}}=-\overline{S}_{2}(1) -
3\overline{S}_{2}(3)
\end{equation}

\noindent and by (\ref{p1}), (\ref{p2}), (\ref{final}), and Proposition \ref{brackets},

\begin{equation} \label{s1}
\overline{S}(1)= g(1) + q^4 y\Sigma(1,0) \frac{(q)_{\infty}
(-q^{25}; q^{25})_{\infty}} {(-q)_{\infty} (q^{25};
q^{25})_{\infty}} - q^{3} \frac{(q^{50}; q^{50})_{\infty}^2
(-q^{10}, -q^{15}; q^{25})_{\infty}} {(q^{10}, q^{40};
q^{50})_{\infty} (-q^{5}, -q^{20}; q^{25})_{\infty}}
\end{equation}

\noindent and

\begin{equation} \label{s3}
\overline{S}(3) = -g(2) + qy^3 \Sigma(2,0) \frac{(q)_{\infty}
(-q^{25}; q^{25})_{\infty}} {(-q)_{\infty} (q^{25};
q^{25})_{\infty}} - q^2 y \frac{(q^{50}; q^{50})_{\infty}^2 (-q^{5},
-q^{20}; q^{25})_{\infty}} {(q^{20}, q^{30}; q^{50})_{\infty}
(-q^{10}, -q^{15}; q^{25})_{\infty}}.
\end{equation}

\noindent By (\ref{gen3}), (\ref{s1}), and (\ref{s3}), we have

$$
\begin{aligned}
& -g(1) - q^4 y\Sigma(1,0) \frac{(q)_{\infty}
(-q^{25}; q^{25})_{\infty}} {(-q)_{\infty} (q^{25};
q^{25})_{\infty}} + q^{3} \frac{(q^{50}; q^{50})_{\infty}^2
(-q^{10}, -q^{15}; q^{25})_{\infty}} {(q^{10}, q^{40};
q^{50})_{\infty} (-q^{5}, -q^{20}; q^{25})_{\infty}} \\
& + 3g(2) - 3qy^3 \Sigma(2,0) \frac{(q)_{\infty}
(-q^{25}; q^{25})_{\infty}} {(-q)_{\infty} (q^{25};
q^{25})_{\infty}} + 3q^2 y \frac{(q^{50}; q^{50})_{\infty}^2 (-q^{5},
-q^{20}; q^{25})_{\infty}} {(q^{20}, q^{30}; q^{50})_{\infty}
(-q^{10}, -q^{15}; q^{25})_{\infty}} \\
& = \Bigl \{ r_{12}(0) q^0 + r_{12}(1)q + r_{12}(2)q^2 +
r_{12}(3)q^3 + r_{12}(4)q^4 \Bigr \}
\frac{(q)_{\infty}}{2(-q)_{\infty}}.
\end{aligned}
$$

\noindent We now multiply the right hand side of the above
expression using \eqref{lem6eq2} and the $R_{12}(d)$ from Theorem
\ref{main}, equating coefficients of powers of $q$.  The
coefficients of $q^{0}$, $q^{1}$, $q^{2}$, $q^{3}$, $q^{4}$ give
us, respectively,

\begin{equation} \label{check0}
\begin{aligned}
3g(2) - g(1) &= 5y^2 \frac{(q^{50}; q^{50})_{\infty}^4 (q^5, q^{10}, q^{40}, q^{45}; q^{50})_{\infty} (q^{25}; q^{25})_{\infty}}{(q^5; q^5)_{\infty}^3 (q^5; q^{10})_{\infty} (-q^{25}; q^{25})_{\infty}} \\
& + y \frac{(q^{25}; q^{25})_{\infty}^2}{(q^5; q^{10})_{\infty}^5 (q^{15}, q^{20}, q^{30}, q^{35}; q^{50})_{\infty} (-q^{25}; q^{25})_{\infty}}
\\ &+ 4y^2 \frac{(q^{50}; q^{50})_{\infty}^3 (q^5, q^{45}; q^{50})_{\infty}}{(q^5; q^{10})_{\infty}^5 (q^{20}, q^{30}; q^{50})_{\infty} (q^{25}; q^{25})_{\infty}} \\
&+ 20y^4 \frac{(q^{50}; q^{50})_{\infty}^7 (q^5, q^{45}; q^{50})_{\infty}^2 (q^{10}, q^{40}; q^{50})_{\infty}^3 (q^5, q^{45}, q^{50}; q^{50})_{\infty}}{(q^5; q^5)_{\infty}^4 (q^{25}; q^{25})_{\infty}^2} \\
&- 4y \frac{(q^{10}; q^{10})_{\infty} (q^{15}, q^{35}, q^{50}; q^{50})_{\infty}}{(q^5; q^{10})_{\infty}^4 (q^{10}, q^{15}, q^{35}, q^{40}; q^{50})_{\infty}^3},
\end{aligned}
\end{equation}

\begin{equation} \label{check1}
\begin{aligned}
\frac{(q^5, q^{10}, q^{40}, q^{45}; q^{50})_{\infty} (q^{15}, q^{35}, q^{50}; q^{50})_{\infty}}{(q^5; q^5)_{\infty} (q^5; q^{10})_{\infty}} &=
y \frac{(q^{50}; q^{50})_{\infty}^2 (q^5, q^{45}; q^{50})_{\infty}^2 (q^{10}, q^{40}; q^{50})_{\infty}^3}{(q^5; q^5)_{\infty}^2} \\
& + \frac{1}{(q^{15}, q^{35}; q^{50})_{\infty}^3 (q^{10}, q^{40}; q^{50})_{\infty}},
\end{aligned}
\end{equation}

\begin{equation} \label{check2}
\begin{aligned}
3 & \frac{(-q^5, -q^{20}; q^{25})_{\infty} (q^{50}; q^{50})_{\infty}^2}{(q^{20}, q^{30} ; q^{50})_{\infty} (-q^{10}, -q^{15}; q^{25})_{\infty}} \\
&=
- 4 \frac{(q^{50}; q^{50})_{\infty}^2 (q^{15}, q^{35}, q^{50}; q^{50})_{\infty}}{(q^5; q^{10})_{\infty}^5 (q^{20}, q^{30}; q^{50})_{\infty} (q^{25}; q^{25})_{\infty}} \\ 
& - 20y^2 \frac{(q^{50}; q^{50})_{\infty}^7 (q^5, q^{45}; q^{50})_{\infty}^2 (q^{10}, q^{40}; q^{50})_{\infty}^3 (q^{15}, q^{35}, q^{50}; q^{50})_{\infty}}{(q^5; q^5)_{\infty}^4 (q^{25}; q^{25})_{\infty}^2} \\ 
&+ 5 \frac{(q^{50}; q^{50})_{\infty}^3 (q^{25}; q^{25})_{\infty}}{(q^5; q^5)_{\infty}^2 (q^{15}, q^{35}; q^{50})_{\infty}^3 (q^5, q^{10}, q^{40}, q^{45}; q^{50})_{\infty} (-q^{25}; q^{25})_{\infty}} \\
&+ 10 \frac{(q^{50}; q^{50})_{\infty}^3 (q^5, q^{45}, q^{50}; q^{50})_{\infty}}{(q^5; q^5)_{\infty}^2 (q^5,q^{45}; q^{50})_{\infty}^3 (q^{15}, q^{20}, q^{30}, q^{35}; q^{50})_{\infty}} \\
&- 8\frac{(q^{25}; q^{25})_{\infty} (q^5, q^{45}, q^{50}; q^{50})_{\infty}}{(q^5; q^{10})_{\infty}^5 (q^{10}, q^{15}, q^{35}, q^{40}; q^{50})_{\infty}},
\end{aligned}
\end{equation}

\begin{equation} \label{check3}
\begin{aligned}
& \frac{(-q^{10}, -q^{15}; q^{25})_{\infty} (q^{50}; q^{50})_{\infty}^2}{(q^{10}, q^{40}; q^{50})_{\infty} (-q^5, -q^{20}; q^{25})_{\infty}} \\
&=
-10y \frac{(q^{50}; q^{50})_{\infty}^3 (q^{15}, q^{35}, q^{50}; q^{50})_{\infty}}{(q^5; q^5)_{\infty}^2 (q^{15}, q^{35}; q^{50})_{\infty}^3 (q^5, q^{10}, q^{40}, q^{45}; q^{50})_{\infty}} \\
& + 5 \frac{(q^{50}; q^{50})_{\infty}^3 (q^{25}; q^{25})_{\infty}}{(q^5; q^5)_{\infty}^2 (q^{5},q^{45}; q^{50})_{\infty}^3 (q^{15}, q^{20}, q^{30}, q^{35}; q^{50})_{\infty} (-q^{25}; q^{25})_{\infty}} \\ 
&- 4\frac{(q^{25}; q^{25})_{\infty}^2}{(q^5; q^{10})_{\infty}^5 (q^{10}, q^{15}, q^{35}, q^{40}; q^{50})_{\infty} (-q^{25}; q^{25})_{\infty}} \\
&+ 4y \frac{(q^{10}; q^{10})_{\infty} (q^5, q^{45}, q^{50}; q^{50})_{\infty}}{(q^5; q^{10})_{\infty}^4 (q^{10}, q^{15}, q^{35}, q^{40}; q^{50})_{\infty}^3},
\end{aligned}
\end{equation}

\begin{equation} \label{check4}
\begin{aligned}
0 &= 5y^2 \frac{(q^{50}; q^{50})_{\infty}^4 (q^5, q^{10}, q^{40}, q^{45}; q^{50})_{\infty} (q^5, q^{45}, q^{50}; q^{50})_{\infty}}{(q^5; q^5)_{\infty}^3 (q^5; q^{10})_{\infty}} \\
&+ y \frac{(q^{25}; q^{25})_{\infty} (q^5, q^{45}, q^{50}; q^{50})_{\infty}}{(q^5; q^{10})_{\infty}^5 (q^{15}, q^{20}, q^{30}, q^{35}; q^{50})_{\infty}} \\ &- 5 \frac{(q^{50}; q^{50})_{\infty}^3 (q^{15}, q^{35}, q^{50}; q^{50})_{\infty}}{(q^5; q^5)_{\infty}^2 (q^{5},q^{45}; q^{50})_{\infty}^3 (q^{15}, q^{20}, q^{30}, q^{35}; q^{50})_{\infty}} \\ 
&+ 4\frac{(q^{25}; q^{25})_{\infty} (q^{50}; q^{50})_{\infty}}{(q^5; q^{10})_{\infty}^5 (q^{10}, q^{40}; q^{50})_{\infty}} \\ &+ \frac{(q^{10}; q^{10})_{\infty} (q^{25}; q^{25})_{\infty}}{(q^5; q^{10})_{\infty}^4 (q^{10}, q^{15}, q^{35}, q^{40}; q^{50})_{\infty}^3 (-q^{25}; q^{25})_{\infty}}.
\end{aligned}
\end{equation}



While we cannot rule out the possibility that \eqref{check0}--\eqref{check4} could be proven with clever applications of identities 
involving infinite products, we verify them using standard computational techniques from the theory of modular forms, which we briefly summarize.  
First, divide each identity by one of its terms to put it in the form
\begin{equation} \label{form}
1 = \sum F_i,
\end{equation}
where each $F_i$ can be expressed in terms of generalized $\eta$-products \cite{Ro1}.  (For \eqref{check0} 
we first use \eqref{g1} and \eqref{g2} to write the left hand side as a sum of $4$ infinite products). By equations (11) and (12) in \cite{Ro1}, one establishes that each of the products $F_i$ is a modular function on some space. (In our case, after letting $q=q^{1/5}$, the functions were always on $\Gamma_1(10)$).  Using Theorem 4 in \cite{Ro1}, we then determine the order of $F_i$ at the cusps and multiply both sides of \eqref{form} by an appropriate power of the Delta function, $\Delta^k(z)$, so that each $\Delta^kF_i$ is holomorphic at the cusps.  Verifying the identity up to $q^T$, where $T$ is the dimension of the appropriate space is then sufficient to confirm its truth.  (In our case, we have
$T = \frac{k(10)^2}{12}\prod_{p \mid 10}(1-1/p^2) = 72k$.)  

We now turn to the rank differences $R_{02}(d)$, proceeding as
above. Again by (\ref{gen1}), (\ref{s}), and (\ref{rels}), we have

\begin{equation} \label{gen4}
\sum_{n=0}^{\infty} \Bigl\{ \overline{N}_{2}(0,5,n) -
\overline{N}_{2}(2,5,n) \Bigr\} q^{n}
\frac{(q)_{\infty}}{2(-q)_{\infty}}=\overline{S}_{2}(5) +
2\overline{S}_{2}(1) + \overline{S}_{2}(3).
\end{equation}

\noindent By Lemma \ref{Sofq} (with $\ell=5$), (\ref{gen4}),
(\ref{s1}), and (\ref{s3}), we have

$$
\begin{aligned}
&  \frac{-(q)_{\infty}}{2(-q)_{\infty}} + \frac{1}{2}  + 2g(1) + 2q^4 y\Sigma(1,0) \frac{(q)_{\infty}
(-q^{25}; q^{25})_{\infty}} {(-q)_{\infty} (q^{25};
q^{25})_{\infty}} - 2q^{3} \frac{(q^{50}; q^{50})_{\infty}^2
(-q^{10}, -q^{15}; q^{25})_{\infty}} {(q^{10}, q^{40};
q^{50})_{\infty} (-q^{5}, -q^{20}; q^{25})_{\infty}}
\\
&  -g(2) + qy^3 \Sigma(2,0) \frac{(q)_{\infty}
(-q^{25}; q^{25})_{\infty}} {(-q)_{\infty} (q^{25};
q^{25})_{\infty}} - q^2 y \frac{(q^{50}; q^{50})_{\infty}^2 (-q^{5},
-q^{20}; q^{25})_{\infty}} {(q^{20}, q^{30}; q^{50})_{\infty}
(-q^{10}, -q^{15}; q^{25})_{\infty}} \\
& = \Bigl \{ r_{02}(0) q^0 + r_{02}(1)q + r_{02}(2)q^2 +
r_{02}(3)q^3 + r_{02}(4)q^4 \Bigr \}
\frac{(q)_{\infty}}{2(-q)_{\infty}}.
\end{aligned}
$$

\noindent Again, equating coefficients of powers of $q$ yields the
following identities.

\begin{equation} \label{check5}
\begin{aligned}
\frac{1}{2} + 2g(1) - g(2) & = \frac{(q^{25};q^{25})_{\infty}^7}{2(q^5;q^{10})_{\infty}^6(q^{15},q^{20},q^{30},q^{35};q^{50})_{\infty}(q^{50};q^{50})_{\infty}^5(-q^{25}; q^{25})_{\infty}} \\ &+
\frac{y^2(q^5;q^5)_{\infty}(q^{50};q^{50})_{\infty}(q^{25}; q^{25})_{\infty}}{2(q^{15},q^{20},q^{30},q^{35};q^{50})_{\infty}^3(q^{25};q^{25})_{\infty}(-q^{25}; q^{25})_{\infty}} \\
&+ \frac{2y(q^{25};q^{25})_{\infty}^3(q^5,q^{45};q^{50})_{\infty}}{(q^5;q^{10})_{\infty}^6(q^{20},q^{30},q^{50};q^{50})_{\infty}} \\ &-
\frac{5y^2(q^{50};q^{50})_{\infty}^3(q^{25}; q^{25})_{\infty}}{(q^5;q^5)_{\infty}^2(q^{10},q^{40};q^{50})_{\infty}(q^{15},q^{35};q^{50})_{\infty}^4(-q^{25}; q^{25})_{\infty}}  \\
&+ \frac{2y^{3} (q^5; q^5)_{\infty} (q^{50}; q^{50})_{\infty}^3 (q^5, q^{45}, q^{50}; q^{50})_{\infty}}{(q^{15}, q^{35}; q^{50})_{\infty}^2 (q^{20}, q^{30}; q^{50})_{\infty}^3 (q^{25}; q^{25})_{\infty}^3} \\
& + \frac{10y (q^{50}; q^{50})_{\infty}^4 (q^{25}; q^{25})_{\infty} (q^{20}, q^{30}; q^{50})_{\infty}}{(q^5, q^{45}; q^{50})_{\infty} (q^5; q^5)_{\infty}^3} \\
& + \frac{2(q^{50}; q^{50})_{\infty} (-q^{10}, -q^{15}; q^{25})_{\infty} (q^{15}, q^{35}, q^{50}; q^{50})_{\infty}}{(-q^5, -q^{20}; q^{25})_{\infty} (q^5, q^{10}, q^{40}, q^{45}; q^{50})_{\infty}} \\
& - \frac{2(q^{25}; q^{25})_{\infty}^2 (q^{10}, q^{40}; q^{50})_{\infty}^5 (q^{50}; q^{50})_{\infty}^2 (q^{15}, q^{35}, q^{50}; q^{50})_{\infty}}{(q^5; q^{10})_{\infty}^4 (q^5, q^{45}; q^{50})_{\infty}^4 (q^{10}; q^{10})_{\infty}^3} \\
& + \frac{2y (q^{25}; q^{25})_{\infty}^4 (q^5, q^{45}; q^{50})_{\infty}}{(q^5; q^{10})_{\infty}^6 (q^{20}, q^{30}; q^{50})_{\infty} (q^{50}; q^{50})_{\infty}^2},
\end{aligned}
\end{equation}

\begin{equation} \label{check6}
0 = 0,
\end{equation}

\begin{equation} \label{check7}
\begin{aligned}
2&\frac{y^2 (q^5; q^5)_{\infty} (q^{50}; q^{50})_{\infty}^4}{(q^{15}, q^{35}; q^{50})_{\infty} (q^{20}, q^{30}; q^{50})_{\infty}^3 (q^{25}; q^{25})_{\infty}^3}  \\
&= \frac{y(q^{50};q^{50})_{\infty}^2(-q^5,-q^{20};q^{25})_{\infty}}{(-q^{10},-q^{15};q^{25})_{\infty}(q^{20},q^{30};q^{50})_{\infty}}  \\
&+ \frac{4y(q^{25}; q^{25})_{\infty} (q^5, q^{45}, q^{50}; q^{50})_{\infty}}{(q^5; q^{10})_{\infty}^5 (q^{10}, q^{15}, q^{35}, q^{40}; q^{50})_{\infty}} \\ &- \frac{5y (q^{50}; q^{50})_{\infty}^3 (q^{25}; q^{25})_{\infty}}{(q^5; q^5)_{\infty}^2 (q^{15}, q^{35}; q^{50})_{\infty}^3 (q^5, q^{10}, q^{40}, q^{45}; q^{50})_{\infty} (-q^{25}; q^{25})_{\infty}},
\end{aligned}
\end{equation}

\begin{equation} \label{check8}
\begin{aligned}
4&\frac{(q^{15}, q^{35}; q^{50})_{\infty}^2 (q^{25}; q^{25})_{\infty}^2}{(q^5; q^{10})_{\infty}^6 (q^{20}, q^{30}; q^{50})_{\infty}} \\
&=  \frac{2(q^{25}; q^{25})_{\infty}^2}{(q^5; q^{10})_{\infty}^5 (q^{10}, q^{15}, q^{35}, q^{40}; q^{50})_{\infty} (-q^{25}; q^{25})_{\infty}} \\
&+ \frac{2 (q^{25}; q^{25})_{\infty}^2 (q^{10}, q^{40}; q^{50})_{\infty}^5 (q^{50}; q^{50})_{\infty}^2 (q^5, q^{45}, q^{50}; q^{50})_{\infty}}{(q^5; q^{10})_{\infty}^4 (q^5, q^{45}; q^{50})_{\infty}^4 (q^{10}; q^{10})_{\infty}^3},
\end{aligned}
\end{equation}

\begin{equation} \label{check9}
\begin{aligned}
&\frac{(q^{25};q^{25})_{\infty}^6 (q^{5}, q^{45}, q^{50}; q^{50})_{\infty}}{(q^5;q^{10})_{\infty}^6(q^{15},q^{20},q^{30},q^{35};q^{50})_{\infty}(q^{50};q^{50})_{\infty}^5} \\
&= 
\frac{-y^2(q^5;q^5)_{\infty}(q^{50};q^{50})_{\infty}(q^{5}, q^{45}, q^{50}; q^{50})_{\infty}}{(q^{15},q^{20},q^{30},q^{35};q^{50})_{\infty}^3(q^{25};q^{25})_{\infty}} \\
&- \frac{4y(q^{25};q^{25})_{\infty}^2(q^5,q^{45};q^{50})_{\infty}(q^{5}, q^{45}, q^{50}; q^{50})_{\infty}}{(q^5;q^{10})_{\infty}^6(q^{20},q^{30},q^{50};q^{50})_{\infty}} \\ &+ 
\frac{10y^2(q^{50};q^{50})_{\infty}^3(q^{5}, q^{45}, q^{50}; q^{50})_{\infty}}{(q^5;q^5)_{\infty}^2(q^{10},q^{40};q^{50})_{\infty}(q^{15},q^{35};q^{50})_{\infty}^4}  \\
&+ \frac{4(q^{25}; q^{25})_{\infty} (q^{15}, q^{35}, q^{50}; q^{50})_{\infty}}{(q^5; q^{10})_{\infty}^5 (q^{10}, q^{15}, q^{35}, q^{40}; q^{50})_{\infty}} \\ &+ \frac{5(q^{50}; q^{50})_{\infty}^3 (q^{25}; q^{25})_{\infty}^2 (q^{20}, q^{30}; q^{50})_{\infty}}{(q^5, q^{45}; q^{50})_{\infty} (q^5; q^5)_{\infty}^3 (q^{15}, q^{35}; q^{50})_{\infty} (-q^{25}; q^{25})_{\infty}} \\
&+ \frac{(q^{50}; q^{50})_{\infty} (-q^{10}, -q^{15}; q^{25})_{\infty} (q^{25}; q^{25})_{\infty}}{y(-q^5, -q^{20}; q^{25})_{\infty} (q^5, q^{10}, q^{40}, q^{45}; q^{50})_{\infty} (-q^{25}; q^{25})_{\infty}} \\ &- \frac{(q^{25}; q^{25})_{\infty}^3 (q^{10}, q^{40}; q^{50})_{\infty}^5 (q^{50}; q^{50})_{\infty}^2}{y(q^5; q^{10})_{\infty}^4 (q^5, q^{45}; q^{50})_{\infty}^4 (q^{10}; q^{10})_{\infty}^3 (-q^{25}; q^{25})_{\infty}}.
\end{aligned}
\end{equation}


These equations were verified using modular forms as with \eqref{check0}--\eqref{check4} above.
\end{proof}

\section{Proof of Corollary \ref{omega}}

\begin{proof}
We first recall the two ``universal mock theta functions" (see Section 6 in \cite{mock})

$$
g_2(x,q):=\frac{(-q)_{\infty}}{(q)_{\infty}} \sum_{n \in \mathbb{Z}} \frac{(-1)^n q^{n(n+1)}}{1-xq^n}
$$

\noindent and

$$
g_3(x,q):=\frac{1}{(q)_{\infty}} \sum_{n \in \mathbb{Z}} \frac{(-1)^n q^{\frac{3}{2}n(n+1)}}{1-xq^n}.
$$

\noindent Watson \cite[p. 66]{Wa1} showed that 

\begin{equation} \label{ww}
\omega(q)=g_3(q, q^2).
\end{equation}

\noindent By Theorem 1.1 in \cite{kang},

\begin{equation} \label{k}
xg_2(x,q)= \frac{\eta^4(2\tau)}{\eta^2(\tau) \vartheta(2\alpha; 2\tau)} + xq^{-1/4} \mu(2\alpha, \tau; 2\tau).
\end{equation}

\noindent Here 

$$
\mu(u, v; \tau):=\frac{a^{1/2}}{\vartheta(v; \tau)} \sum_{n \in \mathbb{Z}} \frac{(-b)^n q^{n(n+1)/2}}{1-aq^n}
$$

\noindent and $\vartheta(v; \tau)$ is the classical theta series with product representation

$$
\displaystyle \vartheta(v; \tau)=q^{1/8} b^{-1/2} \prod_{n=1}^{\infty} (1-q^n)(1-bq^{n-1})(1-b^{-1}q^n),
$$

\noindent where $x=e^{2\pi i\alpha}$, $q=e^{2\pi i \tau}$, $a=e^{2\pi i u}$,  and $b=e^{2\pi i v}$. Now by Theorem 1.1, Theorem 1.1 in \cite{loveoz} and (\ref{k}), we have

\begin{equation} \label{ranks}
\begin{aligned}
& \sum_{n \geq 0} \bigl( \overline{N}_2(0,3,3n+2) - \overline{N}_2(1,3,3n+2) \bigr)q^n - 
\sum_{n \geq 0} \bigl( \overline{N}(0,3,3n+2) - \overline{N}(1,3,3n+2) \bigr)q^n \\
&= \frac{6q(-q^3; q^3)_{\infty}}{(q^3; q^3)_{\infty}} \sum_{n \in \mathbb{Z}} \frac{(-1)^n q^{3n^2 + 6n}}{1-q^{6n+2}} + \frac{6(-q^3; q^3)_{\infty}}{(q^3; q^3)_{\infty}} \sum_{n \in \mathbb{Z}} \frac{(-1)^n q^{3n^2 + 3n}}{1-q^{3n+1}} \\
& = 6q^{-3/4} \mu(2\tau, 3\tau; 6\tau) + 6g_2(q, q^3) \\
& = 12g_2(q, q^3)  - \frac{6\eta^4(6\tau)}{q\eta^2(3\tau) \vartheta(2\tau; 6\tau)}.
\end{aligned}
\end{equation}

\noindent Identity (6.1) in \cite{mock} states that

\begin{equation} \label{id}
g_3(x^4, q^4)= \frac{qg_2(x^6 q, q^6)}{x^2} + \frac{x^2 g_2(x^6 q^{-1}, q^6)}{q} - \frac{x^2 (q^2; q^2)_{\infty}^3 (q^{12}; q^{12})_{\infty} j(x^2 q, q^2) j(x^{12} q^6, q^{12})}{q (q^4; q^4)_{\infty} (q^6; q^6)_{\infty}^2 j(x^4, q^2) j(x^6 q^{-1}, q^2)},
\end{equation}

\noindent where $j(x,q):=(x)_{\infty} (q/x)_{\infty} (q)_{\infty}$. Letting $q \to q^{1/2}$, $x \to q^{1/4}$ in (\ref{id}) and using the fact that $g_2(q^2, q^3)=g_2(q, q^3)$ gives us

\begin{equation} \label{after}
g_3(q, q^2) = 2g_2(q, q^3)  - \frac{(q^6; q^6)_{\infty}^4}{(q^2; q^2)_{\infty} (q^3; q^3)_{\infty}^2}.
\end{equation}

\noindent Thus the result follows after substituting (\ref{after}) into (\ref{ranks}) and using (\ref{ww}).

\end{proof}

\section*{Acknowlegements}
The authors would like to thank Kathrin Bringmann for her comments regarding Lemma \ref{us2}.  The first author was partially
supported by the Agence Nationale de Recherche ANR-08-JCJC-0011.
The second author was partially funded by Science Foundation Ireland 08/RFP/MTH1081.

\end{document}